\documentstyle[art10]{article}
\title{Classification of \\ Plane Congruences in 
 {\bf P}$^4_{\bf C}$ (I)}
\author{Manuel Pedreira 
\and Luis-Eduardo Sol\'{a}-Conde\thanks{Supported by an F.P.U.
fellowship of Spanish Government}}
\date{}
\newtheorem{teorema}{Theorem}[subsection]
\newtheorem{definition}{Definition}[subsection]
\newtheorem{proposition}{Proposition}[subsection]

\newtheorem{lemma}{Lemma}[subsection]

\newtheorem{nota}{Remark}[subsection]

\def\rank{\mathop{\rm rank}}
\def\Hom{\mathop{\rm Hom}}

\begin{document}
\maketitle

{\footnotesize{\bf Authors' address:} Departamento de Algebra, Universidad de Santiago
de Compostela. $15706$ Santiago de Compostela. Galicia. Spain. e-mail: {\tt
pedreira@zmat.usc.es}\\
{\bf Abstract:} Nondegenerate plane congruences in the four--dimensional complex
projective space with degenerate general focal conic are classified by using the focal
method due to Corrado Segre.\\
{\bf Mathematics Subject Classifications (1991):} Primary, 14N05; secondary,
51N35.\\
{\bf Key Words:} Focal points, congruences.}

\vspace{0.1cm}

{\bf Introduction:} The subject of plane congruences in ${\bf P}^4_C$
 is a classical one. An approach to the classification is due to C. Segre 
 by using the so called (algebro--geometric--differential) theory of foci [8].
An attempt to classification of plane congruences in ${\bf P}^4_C$ has
been made by M. Italiani in [3] and [4] where he
basically uses the  more differential geometric tool of E. Cartan's moving
frames. Recently, this subject has been revived its
application to the study of singularities of theta
divisors of Jacobians and to the proof of Torelli's theorem, see [2]. 

\smallskip

We initially
attempted to better understand the classification of plane congruences in ${\bf P}^4_C$
given in [4]. We felt that this classification could be clarified by recovering it
with the focal method, and this is the aim of this paper. Focal methods were introduced
by C. Segre in [7] and [8], and a contemporaneous foundation was given in [1]. The
basic definitions and results about plane congruences in ${\bf P}^4_C$ and their
focal loci can be seen in [2]. A comprehensive exposition of the first
order focal locus appeared in [6].

\smallskip

We study the classification of nondegenerate plane congruences in ${\bf
P}^4_C$ whose focal locus on the general plane consists of a
reducible conic. After proving that the degeneration of the conic is equivalent
to the existence of developable families of planes contained in the congruence and
passing by the plane of the conic (proposition $(1.1.2)$), we show in $1.2$ how the
structure of the focal locus determines the focal conic  and provides the criterium
$1.3$ for the classification. Then, the classification is achieved in $\S 2$ by
doing a very detailed analysis of the first order focal locus. We remark that this is
a birrational classification.

\smallskip

The results on this paper belong to the Ph.D. thesis
of the second author whose advisor is the first one. We thank to Ciro
Ciliberto his criticism of a first version of this paper.

\addtocontents{toc}{\protect\vspace{3ex}}
\section{Preliminaries.}
\addtocontents{toc}{\protect\vspace{1ex}}

Throughout this paper, the base field for algebraic varieties is ${\bf C}$. Let ${\bf
P}^4$ be the $4$--dimensional complex projective space and $G(2,4)$ the
Grassmannian of planes in ${\bf P}^4$. We will consider flat families of planes in
${\bf P}^4$ parametrized by a quasi--projective irreducible algebraic surface
$\Sigma\subset G(2,4)$. They were called plane congruences by classics. We will study
the projective generation of the congruence over the open set $U\subset\Sigma$ where
the first order focal locus has the expected structure. In this sense our
classification is a birrational one and we can suppose that $\Sigma$ is smooth.

\smallskip

Given $\sigma\in\Sigma$, let ${\bf P}^2(\sigma)$ denote the corresponding plane in
${\bf P}^4$, and 
$$
{\cal
I}_{\Sigma}:=\{(P,\sigma)\in {\bf P}^4\times\Sigma;\,P\in {\bf
P}^2(\sigma)\}
$$

is the Incidence Variety. Consider the two projections $p: 
{\cal I}_{\Sigma}\longrightarrow {\bf P}^4$ and $q: 
{\cal I}_{\Sigma}\longrightarrow \Sigma$. We call
$V(\Sigma):=\overline{p({\cal I}_{\Sigma}})$ the projective realization of the
plane congruence. The congruence is said to be nondegenerate if $p$ is
dominant. From now on all congruences will be nondegenerate and we refer to them
with the notation $(\Sigma,p,q)$. 

\smallskip

Recall the definition of focal point and refer to [2] for details: given
$\sigma\in\Sigma$ and a point $Q\in {\bf P}^2(\sigma)$, the differential map
$(dp)_Q:T_{{\cal I}_{\Sigma},(Q,\sigma)}\longrightarrow T_{P^4,p(Q,\sigma)}$
provides the characteristic map 
$$
\chi(\sigma):
T_{\Sigma,\sigma}\longrightarrow H^0({\bf P}^2(\sigma),{\cal
N}_{P^2(\sigma),P^4})
$$
and $Q$ is said to be a focal point of the congruence if there is some $w\in
T_{\Sigma,\sigma}$ such that $\chi(\sigma)(w)(Q)=0$. It's well known how all these
linear maps globalize to a homomorphism $\lambda: T(p)\longrightarrow {\cal N}_{{\cal
I}_{\Sigma},P^4\times\Sigma}$ of sheaves of rank $2$, where
$$
T(p):=\Hom(\Omega_{\Sigma\times P^4|P^4},{\cal 
O}_{\Sigma\times P^4})
$$
and the first order focal locus of the
nondegenerate congruence is  well defined as a proper closed subscheme
$F_1(\Sigma)\subset {\cal I}_{\Sigma}$ by the condition $\det(\lambda)=0$.

\smallskip

We can write the equations of the first order focal locus as it was described in
[6]. First of all, we can compute $\chi(\sigma)(w)(Q)$ by taking a curve
$C\subset\Sigma$ such that $T_{C,\sigma}=\langle w\rangle$ and another curve
$D\subset {\cal I}_{\Sigma}$ passing by $(Q,\sigma)\in {\cal I}_{\Sigma}$ such that
$q:D\longrightarrow C$ is a local isomorphism at $(Q,\sigma)$. If
$T_{D,(Q,\sigma)}=\langle v\rangle$, then $(dq)_{(Q,\sigma)}(v)=w$ and
$\chi(\sigma)(w)(Q)=[(dp)_{(Q,\sigma)}(v)]$, where $[\,]$ denotes the class in
$({\cal N}_{P^2(\sigma),P^4})_Q$.

\smallskip

Let $\Delta\subset{\bf C}$ be a neighbourhood of the origin and $\sigma:
\Delta\longrightarrow C$, $\sigma(0)=\sigma$, a parametrization of $C$ around
$\sigma$. Consider the induced parametrization $Q:\Delta\longrightarrow D$ around
$(Q,\sigma)$, such that $qQ=\sigma$.If $u\in\Delta$ then $Q(u)\in{\bf
P}^2(\sigma(u))$, $(Q,\sigma)=(Q(0),\sigma(0))$, $\sigma'(0)=w$ and $Q'(0)=v$. We
obtain that $Q$ is a focal point for the direction $w\in T_{\Sigma,\sigma}$ iff
$(dp)_{(Q,\sigma)}(v)\in T_{P^2(\sigma),(Q,\sigma)}$, and
this shows that the definition of focal point is justly the one given by C.
Segre [8].

Given a function $G(u,v)$, let $\frac{\partial G}{\partial u}=G_u$. Suppose
$\sigma(u,v)$ a parametrization of $\Sigma$ around $\sigma$,
$\sigma(0,0)=\sigma$. If $w_1=\left.\frac{\partial \sigma}{\partial
u}\right|_{(0,0)}$, $w_2=\left.\frac{\partial \sigma}{\partial
v}\right|_{(0,0)}$, then $T_{\Sigma,\sigma}=\langle w_1,w_2\rangle$. Let
$F_1(u,v)=F_2(u,v)=0$  be the equations of the plane ${\bf P}^2(\sigma(u,v))$ in a
neighbourhood of $(0,0)$. Given a direction $\lambda w_1+\mu w_2$, the plane defined
by the equations 
$$
\lambda F_{1,u}|_{(0,0)}+\mu F_{1,v}|_{(0,0)}=0,\,\,\,\lambda
F_{2,u}|_{(0,0)}+\mu F_{2,v}|_{(0,0)}=0
$$
is called the plane infinitely
near to ${\bf P}^2(\sigma)$ for the direction $\lambda w_1+\mu w_2$, and the focal
points correspond to the intersection of the two planes.

\smallskip

We consider the characteristic map $\chi(\sigma):
T_{\Sigma,\sigma}\longrightarrow H^0({\bf P}^2(\sigma), {\cal N}_{P^2(\sigma),P^4})$.
Since ${\cal N}_{P^2(\sigma),P^4)}\cong ({\cal O}_{P^2(\sigma)}(1))^2$, if 
$\chi(\sigma)(w_1)=(f_{11},f_{12})$ and $\chi(\sigma)(w_2)=(f_{21},f_{22})$, then the
focal points on ${\bf P}^2(\sigma)$ corresponding to the direction $(\lambda:\mu)$ are
the solutions $Q\in{\bf P}^2(\sigma)$ of the
equations $\lambda f_{11}(Q)+\mu f_{21}(Q)=0$, $\lambda f_{12}(Q)+\mu
f_{22}(Q)=0$; i.e. 

$$
\det\left(\begin{array}{cc}

{f_{11}(Q)}&{f_{21}(Q)}\\
{f_{12}(Q)}&{f_{22}(Q)}\\
\end{array}\right)=0
$$

Let $U\subset\Sigma$ be the open set consisting of the smooth points
$\sigma\in\Sigma$ such that $p(F_1(\Sigma))$ does not contain the plane ${\bf
P}^2(\sigma)$. If $\sigma\in U$, the first order focal locus on the plane ${\bf
P}^2(\sigma)$ consists of the conic $f_{11}f_{22}-f_{12}f_{21}=0$. We denote
this conic by $C(\sigma)$, and consider the family $q:F_1(U)\longrightarrow U$.
Let ${\cal H}^1_{2,0}({\bf P}^4)$ denote the Hilbert Scheme of conics in ${\bf
P}^4$. Considering the lower semicontinuous function $U\longrightarrow{\cal
H}^1_{2,0}({\bf P}^4)\longrightarrow{\bf N}$ defined by 
$\sigma\mapsto[C(\sigma)]\mapsto\rank(C(\sigma))$ we obtain the following 

\begin{lemma}
Given a nondegenerate plane congruence in ${\bf P}^4$, either every focal conic is
reducible or the general focal conic is irreducible. \rule{2mm}{2mm}
\end{lemma}  

We will study the plane congruences whose focal conic is degenerate.
These congruences have a geometric characterization that yields a description of 
their projective generation.

\subsection{Developable Families of ${\bf P}^r$'s in ${\bf
P}^n$}\label{Developable}

Let ${\bf P}^n$ be the $n$--dimensional complex projective space and $G(r,n)$ the
Grassmannian parametrizing the $r$--dimensional subspaces  in ${\bf P}^n$. A
$1$--dimensional family of subspaces is said to be a developable family if the
general first order focal locus consists of a ${\bf P}^{r-1}$. In order to
describe which are the $1$--dimensional developable families of planes in ${\bf
P}^n$, consider a smooth curve ${\cal C}\subset G(2,n)$ and its Incidence Variety
${\cal I}_{\cal C}$ with the two projections $p:{\cal I}_{\cal
C}\longrightarrow{\bf P}^n$ and $q:{\cal I}_{\cal C}\longrightarrow{\cal C}$. Let
$U\subset {\cal C}$ be the open set of elements $\sigma\in{\cal C}$ such that
the first order focal locus on the plane ${\bf P}^2(\sigma)$ is the line
$r(\sigma)=\langle Q_1(\sigma),Q_2(\sigma)\rangle$. Write $F_1({\cal
C}):=\overline{\bigcup_{\sigma\in U}r(\sigma)}$ the first order focal scheme.

\begin{proposition}
Given a smooth curve ${\cal C}\subset G(2,n)$ defining a
$1$--dimensional developable family of planes, then this family consists of one of the
next examples:
\begin{enumerate}
\item Planes containing a line, if $\dim p(F_1({\cal C}))=1$.
\item If $\dim p(F_1({\cal C}))=2$, let $F_2({\cal C})$ be the second order
focal scheme of $({\cal C},p,q)$. Then
\begin{enumerate}
\item If $\dim p(F_2({\cal C}))=0$, $\cal C$ consists of the tangent planes to a
cone.
\item If $\dim p(F_2({\cal C}))=1$, $\cal C$ consists of the osculating planes to
a curve.
\end{enumerate}
\end{enumerate}
\end{proposition}

{\bf Proof:} Let $U\subset {\cal C}$ be the open set of elements $\sigma\in
{\cal C}$ such that the first order focal locus of $({\cal C},p,q)$ on the
plane
${\bf P}^2(\sigma)$ is the line $r(\sigma)=\langle
Q_1(\sigma),Q_2(\sigma)\rangle$. We take two curves $D_1,D_2\subset F_1({\cal
C})$ such that $D_i\cap r(\sigma)=Q_i$ and the restriction maps
$q:D_i\longrightarrow{\cal C}$ are local isomorphisms at $Q_i$. Let
$\lambda Q_1+\mu Q_2$ be a point in the line $r(\sigma)$. Since this line is a
focal line for the family of planes, taking a point $Q_3\in{\bf P}^2(\sigma)$ such
that ${\bf P}^2(\sigma)=\langle Q_1,Q_2,Q_3\rangle$, we have  
$\rank(Q_1,Q_2,Q_3,\lambda Q'_1+\mu Q'_2)\leq 3$. So we have two possibilities over
$r(\sigma)$: either $\rank(Q_1,Q_2,\lambda Q'_1+\mu Q'_2)\leq 2$ for every
$(\lambda:\mu)$, or $\rank(Q_1,Q_2,\lambda Q'_1+\mu Q'_2)\leq 2$ for an only one
$(\lambda:\mu)$. In the first case, every point of $F_1({\cal C})$ is a  
focal point, thus the differential map $(dp)_Q: T_{F_1({\cal
C}),(Q,\sigma)}\longrightarrow T_{P^4, p(Q,\sigma)}$ is not injective and $\dim
p(F_1({\cal C}))<2$. So the projective realization of the focal family is a line,
 yielding that our family of planes must be a $1$--dimensional family of planes
containing a line. We remark that this happens just when every  first order focal
point is a second order focal point. For the second case, suppose $P_1\in
r(\sigma)$ is the second order focal point in the general line. Consider the set
$F_2({\cal C})\subset F_1({\cal C})$ of such points. This is a proper closed
subset and we have two possibilities:

\begin{itemize}
\item $p(F_2({\cal C}))$ {\it is a point}. In this case, the projective 
realization of $F_1({\cal C})$ is a cone of vertex $P_1$, and its
generators can be parametrized in the way $\langle
P_1,P_2(u)\rangle$. So the family of planes is of the form $\langle
P_1,P_2(u),P_3(u)\rangle$. Now, the point $P_2$ is focal for this
family, so $\rank\,(P_1,P_2,P_3,P'_2)\,\leq\,3$, but it is
not focal for the family of lines, so
$\rank\,(P_1,P_2,P'_2)\,=\,3$ and ${\bf
P}^2(\sigma)\,=\,\langle P_1,P_2,P'_2\rangle$. That is, 
{\it our family is the family of tangent planes to a cone (consisting of
its foci of first order, and with vertex its focus of the second
order)}.
\item  $p(F_2({\cal C}))$ {\it is a curve}. Since  
$\rank\,(P_1,P_2,P'_1)\,=\,2$, the lines in $F_1({\cal C})$ can
be written as ${\bf P}^1(\sigma(u))\,=\,\langle
P_1(u),P'_1(u)\rangle$. Analogously, one can write the
family of planes ${\bf P}^2(\sigma(u))\,=\,\langle
P_1,P'_1,P''_1\rangle$. So ${\cal C}$ {\it is 
parametrizing the osculating planes to the curve consisting of its foci of the
second order, and ${\cal C}$ is the family of tangent planes to the surface of its
foci of the first order}. \rule{2mm}{2mm}
\end{itemize}

\smallskip

The first step for giving a Criterium to classify the plane congruences in ${\bf
P}^4$ with degenerate focal conic is the next result:

\begin{proposition}
Let $(\Sigma,p,q)$ be a nondegenerate congruence in ${\bf P}^4$. There is a
developable family passing by the general plane ${\bf P}^2(\sigma)$ if and only if
the general focal conic $C(\sigma)\subset{\bf P}^2(\sigma)$ is degenerate.
\end{proposition}

{\bf Proof:} Since the "only if" part is clear, we only prove the "if" part.
Suppose that $C(\sigma)$ is a degenerate conic. If $(\lambda:\mu)\in{\bf P}^1$, let
$C(\sigma,\lambda:\mu)\subset C(\sigma)$ be the set of focal points corresponding to
the direction $(\lambda:\mu)$. If $C(\sigma,\lambda:\mu)$ is a point for all
$(\lambda:\mu)$, we consider the incidence variety ${\cal
I}:=\{(P,(\lambda:\mu))\in {\bf P}^2(\sigma)\times{\bf P}^1;\,P\in
C(\sigma,\lambda:\mu)\}$ with the two projections ${\bar p}: {\cal
I}\longrightarrow{\bf P}^2(\sigma)$ and ${\bar q}: {\cal
I}\longrightarrow{\bf P}^1$, and we see that ${\bar p}({\cal I})=C(\sigma)$ would
be irreducible. Besides, the function ${\bf P}(T_{\Sigma,\sigma})\cong{\bf 
P}^1\longrightarrow{\bf N}$ defined by $(\lambda:\mu)\mapsto\dim{\bar
q}^{-1}(\lambda:\mu)$ is upper semicontinuous and
$$
D_{\sigma}:=\{(\lambda:\mu);\,\dim{\bar q}^{-1}(\lambda:\mu)\geq 1\}
$$
is a closed
subset in ${\bf P}^1$ defined by two quadratic polynomial equations. So, either
$D_{\sigma}$ is a finite subset of degree $1$, $2$; or $D_{\sigma}\cong{\bf P}^1$.
These solutions correspond to developable families passing by the plane ${\bf
P}^2(\sigma)$: they are the curves such that the tangent line at the general point
corresponds to $(\lambda:\mu)\in{\bf P}(T_{\Sigma,\sigma})\cong{\bf P}^1$.
\rule{2mm}{2mm}

\smallskip

{\bf Remark:} The proof of the proposition shows that the developable families
passing by ${\bf P}^2(\sigma)$ correspond to the directions in $T_{\Sigma,\sigma}$
 whose focal locus is a line. Moreover it's clear that there are four
possibilities: there is just one developable family passing by the general
plane; there are two different or coincident developable families passing by the
general plane; and finally, every $1$--dimensional family passing by the general
plane is developable.

\subsection{Projective Generation of the Focal Conic.}\label{generation}

We will show how the focal conic arises in relation to the directions
$(\lambda:\mu)$. We observe that the focal conic is generated by the two pencils of
lines in ${\bf P}^2(\sigma)$ corresponding to the traces of the two pencils of
hyperplanes $\lambda f_{11}+\mu f_{21}=0$ and $\lambda f_{12}+\mu f_{22}=0$. We have
the following possibilities:

\begin{itemize}
\item {\bf Only one of the pencils degenerates on a line not belonging to the other
pencil}: In this case, for example, there is some $a\neq 0$ such that
$f_{11}=af_{21}$;
$f_{12}$ and $f_{22}$ are not proportional and $f_{21}$ is not a linear
combination of $f_{12}$ and $f_{22}$. The focal conic is 
$f_{21}(af_{12}-f_{22})=0$ where $af_{12}-f_{22}=0$ is the focal line
corresponding to the direction $(a:-1)$ and for other direction
$(\lambda:\mu)\neq (a:-1)$, the focal locus consists of the point
$\lambda f_{12}+\mu f_{22}=f_{21}=0$. So, we only have one developable system
and the focal loci are disjoint.
\item {\bf The two pencils degenerate on two lines}: In this case there is some $a\neq
0$ and some $b\neq 0$ such that $f_{11}=af_{21}$, $f_{12}=bf_{22}$, and the focal
conic is $f_{21}f_{22}=0$. We can suppose $a\neq b$ since, in other case, every
point of the plane would be focal for $(\lambda:\mu)=(a:-1)$. First, we suppose
$f_{21}$ and $f_{22}$ are not proportional. We get the focal line $f_{22}=0$
corresponding to $(\lambda:\mu)=(a:-1)$ and the focal line $f_{21}=0$ 
corresponding to $(\lambda:\mu)=(b:-1)$. For other values of $(\lambda:\mu)$, the
focal locus is the intersection point $f_{21}=f_{22}=0$. So we have two different
developable systems and a focal point for every direction. Finally, if
$f_{21}$ and $f_{22}$ are proportional, the focal conic consists of the double line
$(f_{21})^2=0$, and every point is focal for every direction.
\item {\bf Two pencils with different base points}: We call $Q_1$ and $Q_2$ the base
points of the two pencils $H_{Q_1}$ and $H_{Q_2}$. Since the conic is degenerate,
 there is a value $(\lambda:\mu)$ such that the equations $\lambda f_{11}+\mu
f_{21}=0$ and $\lambda f_{12}+\mu f_{22}=0$ give the same line, $\langle Q_1,
Q_2\rangle$. Suppose that this value is $(1:0)$, there is some $a\neq 0$ such that
$f_{12}=af_{11}$ and the focal conic has equation $f_{11}(f_{22}-af_{21})=0$. We get
the focal line $f_{11}=0$ corresponding to $(\lambda:\mu)=(1:0)$, and if
$(\lambda:\mu)\neq(1:0)$, the focal locus is the solution of $\lambda f_{11}+\mu
f_{21}=f_{22}-af_{21}=0$. The focal points are on the line $f_{22}-af_{21}=0$ but not
on $f_{11}=0$. This case is analogous to the first one: there is only one developable
system passing by the general plane and the focal loci are disjoint. 
\item {\bf Two pencils with the same base point, and one of the pencils can be
degenerate}. Suppose that the first pencil $H_1$ is nondegenerate. If
$f_{12}=\lambda_1f_{11}+\mu_1f_{21}$, $f_{22}=\lambda_2f_{11}+\mu_2f_{21}$, the
focal conic is defined by the equation $\lambda_2
f_{11}^2+(\mu_2-\lambda_1)f_{11}f_{21}-\mu_1f_{21}^2=0$. The base point is focal
for every direction $(\lambda:\mu)$ except for 
$(\lambda:\mu)=(\lambda\lambda_1+\mu\lambda_2:\lambda\mu_1+\mu\mu_2)$. That's to
say, except for the eigenvectors of the matrix $A=\left(
^{\lambda_1\,\lambda_2}_{\mu_1\,\mu_2}\right)$. We have the following possibilities:
\begin{itemize}
\item[-] If $H_2$ is degenerate, then $A$ is degenerate and, choosing a value
$(\lambda_0:\mu_0)\in {\bf P}^1$ such that $\lambda_0\lambda_1+\mu_0\lambda_2=
\lambda_0\mu_1+\mu_0\mu_2=0$, we get the line $\lambda_0f_{11}+\mu_0f_{21}=0$.
Taking other different value $(\lambda':\mu')$, we get the line
$\lambda'f_{11}+\mu'f_{21}=0$ and the focal conic is
$(\lambda_0f_{11}+\mu_0f_{21})(\lambda'f_{11}+\mu'f_{21})=0$.
\item[-] If $H_2$ is nondegenerate, we have two possibilities: either both pencils
have common two different values $(\lambda:\mu)$, the eigenvectors of $A$, and the
focal conic is as above; or they have in common all the values, and thus 
$\mu_1=\lambda_2=0$ and $\lambda_1=\mu_2$ and every point of the plane is a focal
point.
\item[-] Finally, it is possible that, being $H_2$ nondegenerate, the matrix
$A$ has a unique eigenvector $(\lambda_0:\mu_0)$. Then the focal conic is
the double line $(\lambda_0 f_{11}+\mu_0 f_{21})^2=0$. In this case
$(\lambda_0:\mu_0)$ is a double developable direction and, for the other
values, the focal locus consists of the point $f_{11}=f_{21}=0$. Now, the existence
of this focal conic is equivalent to the condition
$(\mu_2-\lambda_1)^2+4\mu_1\lambda_2=0$.
\end{itemize}
\end{itemize}

\subsection{Criterium for the Classification.}\label{Crite}

If $(\Sigma, p, q)$ is a nondegenerate plane congruence such that the
general focal locus is a degenerate conic, the congruence belongs to one of the
following types:\\
$\alpha$--{\bf Congruence:} There is only one developable system passing through 
the general plane. Equivalently, the general
focal conic consists of two different lines $r\vee r'$ and the focal
loci corresponding to all directions are disjoint. For only one direction we get
the focal line $r$, and for other different one we get a point in
$r'\setminus(r\cap r')$.\\
$\beta$--{\bf Congruence:} There are two different developable systems passing 
through the general plane. Equivalently, the general
focal conic consists of two different lines $r\vee r'$ with $r$ and
$r'$ corresponding to two different developable directions. For other direction
different from these ones, the focal locus consists of the singular point $r\cap r'$
which is a focal point for every direction.\\
$\gamma$--{\bf Congruence:} There are two coincident developable systems
passing through the general plane. Equivalently, the
general focal conic consists of a double line $r^2$, with $r$ the focal locus for
a double developable direction and for all other directions, the focal locus is a
fixed point $P$, which is a focal point for all of them.\\
$\delta$--{\bf Congruence:} Every $1$--dimensional family passing through the general
plane of the congruence is developable. This is equivalent to the property that the
general focal conic is a double line
$r^2$, being $r$ a focal line for every direction.

\addtocontents{toc}{\protect\vspace{3ex}}
\section{The Classification.}
\addtocontents{toc}{\protect\vspace{1ex}}

\subsection{$\delta$--Congruences}\label{delta}
Let $(\Sigma,p,q)$  be a $\delta$--congruence and consider the 
open set $U\subset\Sigma$ of $\sigma\in U$ satisfying $C(\sigma)=r(\sigma)^2$, where 
$r(\sigma)$ focal line for every direction.

\begin{teorema}
If $(\Sigma,p,q)$ is a $\delta$--congruence, it consists of the linear system of
planes containing a line $r\subset{\bf P}^4$. This line is the focal locus of the
congruence on every plane.
\end{teorema}

{\bf Proof:} We define a morphism $\psi: U\longrightarrow G(1,4)$ by
$\psi(\sigma)=[r(\sigma)]$ and we will prove that $\Sigma':=\overline{\psi(U)}$ is of
 dimension $0$. Suppose $\dim\Sigma'=2$. If $r(\sigma)\subset{\bf P}^2(\sigma)$
is a focal line and $x\in r(\sigma)$ a general point, we take parametrizations
$x,y,z:\Delta\subset{\bf C}^2\longrightarrow{\cal I}_U$ and
$\sigma:\Delta\longrightarrow U$ such that $\sigma(0,0)=\sigma$,
$x(0,0)=x$, $r(\sigma(u,v))=\langle x(u,v),y(u,v)\rangle$, and ${\bf
P}^2(\sigma(u,v))=\langle x(u,v),y(u,v),z(u,v)\rangle$. Since $x$ is a focal point
for every direction, $\rank(x,y,z,\lambda x_u+\mu x_v)=3$ for every $(\lambda:\mu)$.
So there is $(\lambda_0:\mu_0)$ such that $x$ is a focal point for the family of
lines $\{r(\sigma)\}_{\sigma\in U}$ for the corresponding direction. But, in
particular, this implies that
$\dim p(F_1(U))<\dim F_1(U)$ and $p(F_1(U))$ must be a plane. This is false since
every developable family parametrized by a curve in $U$ would consist of the tangent
planes to a $1$--dimensional family of lines in the plane. Now, suppose that
$\dim\Sigma'=1$. In this case $p(F_1(U))$ is a ruled surface and there is an
infinity of planes in the congruence passing through every generator. We can take a
curve $\gamma\subset U$ such that $\overline{\psi(\gamma)}=\Sigma'$. Since the focal
loci of the planes parametrized by $\gamma$ consists of the generators of the ruled
surface $p(F_1(U))$, by proposition $(1.1.1)$, $\gamma$ will be the family of tangent
planes to $p(F_1(U))$, and so $\dim U=1$. This is false too. \rule{2mm}{2mm}

\subsection{$\beta$/$\gamma$--Congruences}\label{bega}
Let $(\Sigma,p,q)$  be a $\beta$/$\gamma$--congruence. Recall that, if
$\sigma\in\Sigma$ is a general point, there is only one point $P(\sigma)\in
C(\sigma)$ that is a focal point for every direction. Equivalently,
$\rank(dp)_{P(\sigma)}=2$. Consider the closed subset $S(U):=\{(P,\sigma)\in{\cal
I}_U:\,\rank(dp)_{P(\sigma)}\leq 2\}\subset {\cal I}_U$. We can suppose, 
 by restricting $U$ if necessary, that $\rho^{-1}=q|_{S(U)}\longrightarrow U$
is an isomorphism. Consider $\phi:= p\rho$, $\phi(\sigma)=P(\sigma)$, and let 
$\Sigma':=\overline{\phi(U)}$. For a general point $\sigma\in U$,
$(d\phi)_{\sigma}$ is surjective, yielding $T_{\Sigma',P(\sigma)}\subset
T_{P^2(\sigma),P(\sigma)}$. There are three possibilities in according to the
dimension of $\Sigma'$.

\smallskip

{\bf $\beta_1$/$\gamma_1$--Congruences: $\dim\Sigma'=2$.} In this case, the
congruence $(\Sigma,p,q)$ consists of the tangent planes to the surface $\Sigma'$.
We study how the intrinsic geometry of $\Sigma'$ provides a $\beta$ or
$\gamma$--congruence.

\begin{definition} A curve $C$ on a surface $S$ is said to be an Asymptotic Curve,
when the osculating plane to $C$ at $P$, $T_{2,C,P}$, and the tangent plane to $S$
at $P$, $T_{S,P}$, are coincident for every point of $C$. Two $1$-dimensional
families of curves $\gamma$ and $\gamma'$ on a surface $S$ are said to be a
Conjugate Double System if the ruled surface of the tangent lines to
the curves of the family $\gamma'$ at the points of a curve $C\in\gamma$ is
developable; that is, the general focal locus is a point.
\end{definition}

\begin{nota} If $x(u,v)$ denotes a parametrization of $S$ around a general
point, it's well known that (see [5]):
\begin{enumerate}
\item there is a $1$--dimensional family of asymptotic curves on $S$ with
differential equation $dv-\lambda du=0$ if and only if there is a function
$\lambda(u,v)$ such that $\rank(x,x_u,x_v,x_{uu}+(2\lambda)x_{uv}+\lambda^2
x_{vv})=3$ for every $(u,v)$,
\item there is a conjugate double system on $S$ with differential equation
$(dv-\lambda du)(dv-\mu du)=0$ if and only if there are functions $\lambda(u,v)$
and $\mu(u,v)$ such that $\rank(x,x_u,x_v,x_{uu}+(\lambda+\mu)x_{uv}+\lambda\mu
x_{vv})=3$ for every $(u,v)$.
\end{enumerate}
\end{nota}

\begin{teorema} A $\beta_1$/$\gamma_1$--congruence $(\Sigma,p,q)$ consists of the
tangent planes to a nondegenerate irreducible surface $\Sigma'\subset{\bf P}^4$,
and either $\Sigma'$ has a conjugate double system and, equivalently, $(\Sigma,p,q)$
is a $\beta_1$--congruence, or $\Sigma'$ only has a $1$--dimensional family of
asymptotic curves and $(\Sigma,p,q)$
is a $\gamma_1$--congruence.
\end{teorema}

{\bf Proof:} We lift the parametrization $x(u,v)$ of $\Sigma'$ to a
parametrization $\sigma:\Delta\subset{\bf C}^2\longrightarrow U$ such that ${\bf
P}^2(\sigma(u,v))=\langle x(u,v),x_u(u,v),x_v(u,v)\rangle$ is the tangent plane
to $\Sigma'$ at the point $P(\sigma)=x(u,v)$. The condition in Remark $(2.2.1.1)$ is
equivalent to
$(\frac{\partial}{\partial u}+\lambda\frac{\partial}{\partial v})(x_u+\lambda
x_v)\in{\bf P}^2(\sigma(u,v))$, and $\Sigma'$ has a $1$--dimensional family of
asymptotic curves iff $x_u+\lambda x_v$ is a focal point for
$(\frac{\partial}{\partial u}+\lambda\frac{\partial}{\partial v})$ on every
plane. Analogously, $\Sigma'$ has a conjugate double system iff $x_u+\lambda x_v$ is
a focal point for $(\frac{\partial}{\partial u}+\mu\frac{\partial}{\partial v})$
and, symmetrically, $x_u+\mu x_v$ is
a focal point for $(\frac{\partial}{\partial u}+\lambda\frac{\partial}{\partial
v})$. Since $P(\sigma)=x$ is a focal point for every direction, the focal conic
consists of $\langle x,x_u+\lambda x_v\rangle\vee\langle x,x_u+\mu x_v\rangle$,
and $(\Sigma,p,q)$ will be a $\beta_1$--congruence.\\ 
If $\Sigma'$ has a
$1$--dimensional family of asymptotic curves, $\langle x,x_u+\lambda x_v\rangle$
is a focal line on every plane. If $(\Sigma,p,q)$ is a $\beta_1$--congruence,
there is a point $x_u+\lambda' x_v$, $\lambda'\neq\lambda$ which is a focal
point for $(\frac{\partial}{\partial u}+\mu'\frac{\partial}{\partial v})$. By the
condition in Remark $(2.2.1.2)$, there would be a third focal point  $x_u+\mu' x_v$,
then $\mu'=\lambda'$. So the differential equation $dv-\lambda' du=0$ provides another
$1$--dimensional family of asymptotic curves on $\Sigma'$. We prove that if
$\Sigma'\subset{\bf P}^4$ is a surface with two families of asymptotic curves, the
congruence $\Sigma$ of tangent planes to $\Sigma'$ is degenerate. In fact, for the
general point $\sigma\in\Sigma$, every point in ${\bf P}^2(\sigma)$ is a focal point.
Suppose a parametrization of $\Sigma'$ such that the asymptotic curves have equations
$dv=0$, $du=0$, then $x_{uu},x_{vv}\in\langle x,x_u,x_v\rangle$ for every $(u,v)$ and
every point $Bx_u+Cx_v\in\langle x_u,x_v\rangle$ is a focal point for the direction
$(-B\frac{\partial}{\partial u}+C\frac{\partial}{\partial v})$. \rule{2mm}{2mm}

\smallskip

{\bf $\beta_2$/$\gamma_2$--Congruences: $\dim\Sigma'=1$.} In this case,
$\dim\phi^{-1}(P(\sigma))=1$ for the general $P(\sigma)\in\Sigma'$ and the
congruence $(\Sigma,p,q)$ consists of $1$--dimensional families of planes passing 
through the tangent lines to $\Sigma'$. 

\begin{teorema} The congruence $(\Sigma,p,q)$ is a
$\beta_2$/$\gamma_2$--congruence iff it consists of $1$--dimensional families of
planes passing through the tangent lines to a curve $\Sigma'$. $(\Sigma,p,q)$ is a
$\gamma_2$--congruence iff all the families containing the tangent line are linear
ones and they contain the corresponding osculating plane to $\Sigma'$.
\end{teorema}

{\bf Proof:} We observe that there is a developable family
$\gamma_{P(\sigma)}=\{\tau\in U:\,\phi(\tau)=P(\sigma)\}$ passing through the plane
${\bf P}^2(\sigma)$. We study the existence of other different or coincident
developable system passing through ${\bf P}^2(\sigma)$. The idea consists of taking a
plane $\pi$ that does not meet the curve $\Sigma'$ and analyzing the intersection
with the general plane of $(\Sigma,p,q)$ in order to be able to parametrize the
congruence by the points of $\pi$. Considering a general plane $\pi\subset{\bf
P}^4$, there is an open subset $U\subset\Sigma$ such that ${\bf P}^2(\sigma)\cap\pi$
is a point for $\sigma\in U$. We have a domi\-nant morphism
$\psi:U\longrightarrow\pi$, $\psi(\sigma)={\bf P}^2(\sigma)\cap\pi$ and ${\bf
P}^2(\sigma)=\langle T_{\Sigma',\phi(\sigma)},\psi(\sigma)\rangle$. Since $\psi$
is a local isomorphism at the general point $\sigma\in U$, we can take a
parametrization $\sigma(u,v)$ of $U$ around $\sigma$ and a reference system on
$\pi$ such that $\sigma(0,0)=\sigma$; $\psi(\sigma(u,v))=(u:v:1)$. Let 
$C_{P(\sigma)}=\psi(\gamma_{P(\sigma)})$ the curves satisfying the
differential equation $dv-kdu=0$, for some function $k(u,v)$, around the point
$\psi(\sigma)$. The tangent vectors at every point are $\frac{\partial}{\partial
u}+k\frac{\partial}{\partial v}$ and
$(d\phi)_{\sigma}(\frac{\partial\sigma}{\partial u}+k\frac{\partial\sigma}{\partial
v})=0$, then we complete it with $\frac{\partial}{\partial v}$ to get a base of
$T_{\pi,\psi(\sigma)}$, and ${\bf
P}^2(\sigma(u,v))=\langle\phi(\sigma(u,v)),\frac{\partial(\phi\sigma)}{\partial
v},\psi(\sigma(u,v))\rangle$. There is a developable direction different from $\frac{\partial}{\partial
u}+k\frac{\partial}{\partial v}$ in $T_{\Sigma,\sigma}$ iff there are $\lambda$
and $(a:b)\in{\bf P}^1$ such that:
$\rank(\phi(\sigma),\frac{\partial(\phi\sigma)}{\partial
v},\psi(\sigma),(\lambda(\frac{\partial}{\partial u}+k\frac{\partial}{\partial
v})+\frac{\partial}{\partial v})(a\frac{\partial (\phi\sigma)}{\partial
v}+b\psi(\sigma)))\leq 3$. Equivalently, there is
$\lambda$ such that
$$
\rank(\psi(\sigma),\phi(\sigma),\frac{\partial(\phi\sigma)}{\partial
v},\frac{\partial^2(\phi\sigma)}{\partial
v^2},\lambda(\frac{\partial(\psi\sigma)}{\partial
u}+k\frac{\partial(\psi\sigma)}{\partial v})+\frac{\partial(\psi\sigma)}{\partial
v})\leq 4
$$ 
That is: there is a line
$\langle\psi(\sigma),\lambda(\frac{\partial(\psi\sigma)}{\partial
u}+k\frac{\partial(\psi\sigma)}{\partial v})+\frac{\partial(\psi\sigma)}{\partial
v}\rangle\subset\pi$ different from $T_{C_{P(\sigma)},\psi(\sigma)}$ meeting
$T_{2,\Sigma',P(\sigma)}=\langle\phi(\sigma),\frac{\partial(\phi\sigma)}{\partial
v},\frac{\partial^2(\phi\sigma)}{\partial
v^2}\rangle$. Observe that $T_{2,\Sigma',P(\sigma)}$ meets  $\pi$ in a
point $R(P(\sigma))$. So there will not be other developable direction (and
$(\Sigma,p,q)$ will be a $\gamma_2$--congruence) if and only if $R(P(\sigma))$
belongs to every tangent line to $C_{P(\sigma)}$, and thus these curves are lines
passing through $R(P(\sigma))=\pi\cap T_{2,\Sigma',P(\sigma)}$. \rule{2mm}{2mm}

\smallskip

{\bf $\beta_3$/$\gamma_3$--Congruences: $\dim\Sigma'=0$.} In this case all 
planes of the congru\-ence meet the point $\Sigma'$, and we get the next
characterization:

\begin{teorema} A $\beta_3$/$\gamma_3$--congruence consists of the cone, with
vertex a point $\Sigma'$, over a line congruence $\overline{\Sigma}$ in a hyperplane
$H\not\ni\Sigma'$. $(\Sigma,p,q)$ is $\beta_3$ iff $\overline{\Sigma}$ is not
parabolic; that is, the focal locus on the general line consists of two
points. $(\Sigma,p,q)$ is $\gamma_3$ iff $\overline{\Sigma}$ is
parabolic; that is, the focal locus on the general line consists of one
point.
\end{teorema}

{\bf Proof:} Let $H\subset{\bf P}^4$ be a hyperplane such that $\Sigma'\not\in H$.
We get an isomorphism $\psi:\Sigma\longrightarrow\overline{\Sigma}\subset G({\bf
P}^1,H)$, $\sigma\mapsto[{\bf P}^2(\sigma)\cap H]=[{\bf P}^1(\sigma)]$. To conclude, 
it's enough to observe that $Q\in{\bf P}^1(\sigma)$ is a focal point for the line
congruence iff $\langle Q,P(\sigma)\rangle$ is a focal line for $\Sigma$.
\rule{2mm}{2mm}

\subsection{$\alpha$--Congruences}\label{alfa}
Let $U\subset\Sigma$ be the open set such that the focal locus on ${\bf
P}^2(\sigma)$, $\sigma\in U$, consists of the conic $r(\sigma)\vee r'(\sigma)$.
Suppose $r(\sigma)$ is the focal line corresponding to the only developable system
$\gamma$ passing through ${\bf P}^2(\sigma)$. We stratify $U=\bigcup_{i\in
I}\gamma_i$, where  $\gamma_i$ are the developable systems contained in the
congruence. Let $R(U):=\{(P,\sigma):\,P\in r(\sigma)\}=\bigcup
F_1(\gamma_i)\subset F_1(U)$ and consider the restrictions ${\bar
p}:R(U)\longrightarrow{\bf P}^4$,
${\bar q}:R(U)\longrightarrow U$ of the maps $p$, $q$.
\begin{enumerate}
\item If $\dim{\bar p}(R(U))=2$, we get that $p(F_1(\gamma_i))$ must be
$1$--dimensional by arguing as in theorem $2.1.1$. So ${\bar p}(R(U))$
is a ruled surface with generators  $p(F_1(\gamma_i))$.
\item If $\dim{\bar p}(R(U))=3$, $p(F_1(\gamma_i))$ must be $2$--dimensional for
the general curve $\gamma_i$. So $p(F_1(\gamma_i))$ is a developable ruled
surface and $\gamma_i$ is the corresponding family of tangent planes. In
particular, $p(F_1(\gamma_i))$ is a cone if $p(F_2(\gamma_i))$ is a point. Now let
$(P,\sigma)\in r(\sigma)\subset F_1(\gamma_i)\subset R(U)$ be a general point.
Since $r(\sigma)$ is the focal line in ${\bf P}^2(\sigma)$ for the direction
$T_{\gamma_i,\sigma}$, we get $\ker(dp)_{(P,\sigma)}|_{T_{R(U),(P,\sigma)}}\subset
T{\gamma_i,\sigma}$. Then, because $(d{\bar
p})_{(P,\sigma)}=(dp)_{(P,\sigma)}|_{T_{R(U),(P,\sigma)}}$, $(P,\sigma)$ is focal
 in $R(U))$ iff $(P,\sigma)$ is focal in $\gamma_i$, and $F_1(R(U))=\bigcup
F_2(\gamma_i)$; that is, the focal locus of $R(U)$ consists of a point over
each general line. Since $(\Sigma,p,q)$ is not a $\beta_3$ or
$\gamma_3$--congruence, necessarily $\dim p(F_1(R(U)))>0$. We have two
possibilities: \\either
$\dim p(F_1(R(U)))=2$ and $\dim p(F_1(\gamma_i))=1$ ($\gamma_i$ consists of the
osculating planes to $p(F_2(\gamma_i))$), or
$\dim p(F_1(R(U)))=1$ and
$\dim p(F_1(\gamma_i))=0$ ($\gamma_i$ consists of the
tangent planes to the cone $p(F_1(\gamma_i))$). 
\end{enumerate}

We have the following cases 

\smallskip

{\bf $\alpha_1$--Congruences:} $\dim p(R(U))=3$ and $\dim p(F_1(R(U)))=2$.

\begin{teorema} $(\Sigma,p,q)$ is an $\alpha_1$--congruence iff it consists of the
family of osculating planes to a $1$--dimensional family of curves $C_i$ on a
nondegenerate and irreducible surface $S\subset{\bf P}^4$, except if $C_i$ is
asymptotic ($\Sigma$ is a $\gamma_1$--congruence) or $\{C_i\}$ is one of the
families of a conjugate double system ($\Sigma$ is a $\beta_1$--congruence).
\end{teorema}

{\bf Proof:} Given an $\alpha_2$--congruence, the result follows by taking
$C_i=p(F_2(\gamma_i))$ on the nondegenerate and irreducible surface
$S=p(F_1(R(U)))$. Conversely, given $S$ and $C_i$ on $S$, let $\Sigma$ be the
congruence defined by the osculating planes. If $\Sigma$ is not an
$\alpha$--congruence, there is a point in each plane which is  a focal point for
every direction. Suppose a general point $x\in S$ and let $x(u,v)$ be a
parametrization of $S$ around $x$ such that the curves $C_i$ are the solutions of
$dv=0$. We get a parametrization $\sigma(u,v)$ of $\Sigma$ such that ${\bf
P}^2(\sigma(u,v))=\langle x(u,v),x_u(u,v),x_v(u,v)\rangle$. $\langle x,x_u\rangle$
is focal for the direction $\frac{\partial}{\partial u}$. Assume the
existence of a point in $\langle x,x_u\rangle$ that is a focal point for
$\frac{\partial}{\partial v}$. This is equivalent to the condition
$\rank(x,x_u,x_{uu},x_v,x_{uv})\leq 4$; that is, there is $(c:d)$ such that
$\rank(x,x_u,x_v,cx_{uu}+dx_{uv})\leq 3$ for every $(u,v)$. This provides, by Remark
$2.2.1$ the exceptions in the theorem. \rule{2mm}{2mm}

\smallskip

{\bf $\alpha_2$--Congruences:} $\dim p(R(U))=3$ and $\dim p(F_1(R(U)))=1$.

\begin{teorema} $(\Sigma,p,q)$ is an $\alpha_2$--congruence iff it consists of the
tangent planes to a $1$--dimensional family of cones whose vertexes vary along a
curve $C\subset{\bf P}^4$, and such that the focal locus of the generators of
each cone consists of a point on every line.
\end{teorema}

{\bf Proof:} If $(\Sigma,p,q)$ is an $\alpha_2$--congruence, we conclude by taking
the cones $Q_i:=p(F_1(\gamma_i))$ whose vertexes vary  along the curve
$C:=p(F_1(R(U)))$ and the focal locus of the generators of $Q_i$ is a point on
each generator. Conversely, given a curve $C$ and a $1$--dimensional family of
cones $Q_i$, whose vertexes vary along C, we consider the family $\Sigma$ of their
tangent planes. Since the tangent planes to each cone is a developable
system, the congruence is not in cases $\delta$, $\gamma_1$ and $\gamma_2$; and since
the vertexes of the cones vary along a curve, $\Sigma$ cannot be a $\beta_3$ or
$\gamma_3$ congruence. Then the possibilities are: $\alpha$, $\beta_1$ or
$\beta_2$.

\begin{itemize}
\item If $\Sigma$ is a $\beta_1$--congruence, it consists of the tangent planes to
a surface $\Sigma'$ with a conjugate double system. Consider a parametrization
such that the conjugate double system is defined by the equation $(dv)(du)=0$. The
family of lines $R(U)$ is parametrized by $\langle x(u,v),x_u(u,v)\rangle$ and
there are two different focal points on each line: the point $x$ for the
direction $\frac{\partial}{\partial u}$; and the vertex of the corresponding cone 
$ax+bx_u$, $b\neq 0$, for the direction $\frac{\partial}{\partial v}$.
\item If $\Sigma$ is a $\beta_2$--congruence, there is another different developable
system passing through each plane ${\bf P}^2(\sigma)$: it consists of planes
containing the line $T_{\Sigma',P(\sigma)}$ meeting to $r(\sigma)$ in $P(\sigma)$. We
see that $\Sigma'$ is a directrix for each cone and all of its points are
fundamental points, then they are focal points for the family $\{r(\sigma)\}$.
 Furthermore the vertex of the cone is a focal point on every line, and there are two
focal points on every line. 
\end{itemize}

We conclude that $\Sigma$ is an $\alpha_2$--congruence. \rule{2mm}{2mm}

\smallskip

{\bf $\alpha_3$--Congruences:} $\dim p(R(U))=2$.

\begin{teorema} $(\Sigma,p,q)$ is an $\alpha_3$--congruence iff it consists of a
$1$--dimensional family of planes passing through every generator of a
nondegenerate and irreducible ruled surface $S$ except the tangent planes to
$S$, which is a
$\gamma_1$--congruence.
\end{teorema}

{\bf Proof:} If $(\Sigma,p,q)$ is a $\alpha_3$--congruence,
$S:=\overline{p(R(U))}$ is a ruled surface with generators $r_i:=p(F_1(\gamma_i))$
and there is a $1$--dimensional family of planes passing through every
generator.
$S$ cannot be a developable surface because $\Sigma$ would be in cases
$\beta_2$,
$\gamma_2$, $\beta_3$ or $\gamma_3$. Conversely, given such a surface
$S\subset{\bf P}^4$ and the congruence $\Sigma$ constructed in that
way, is not in case $\delta$. Suppose $\Sigma$ is not an
$\alpha$--congruence and consider
$\Sigma':=\overline{p(S(U))}\subset\overline{p(R(U))}=S$. Since $S$ is neither
the tangent developable to a curve nor a cone, $\Sigma$ cannot be in cases
$\beta_2$/$\gamma_2$ and $\beta_3$/$\gamma_3$. So $\Sigma'=S$ and
$\Sigma$ is the family of tangent planes to $S$. \rule{2mm}{2mm}

\addtocontents{toc}{\protect\vspace{3ex}}
\section{Conclusions.}
\addtocontents{toc}{\protect\vspace{1ex}}

\begin{tabular}{|p{2.7cm}|p{0.3cm}|p{7.7cm}|}
\hline
{\footnotesize \bf Criterium $1.3$}& { }& {\footnotesize \bf Characterization}\\
\hline
{ \quad}& {\footnotesize $\alpha_1$}& {\footnotesize Family of
osculating planes to a $1$--dimensional family
of curves on an irreducible nondegenerate surface $S$
($=p(F_1(R(U)))$), except the case where such curves are
asymptotic or they are curves of a conjugate double system
of $S$.}\\ \cline {2-3}
 {\footnotesize $\alpha$: Only one developable system passing through the
general plane.}&
 {\footnotesize $\alpha_2$}& {\footnotesize Family of tangent planes
to a $1$-dimensional family of cones with vertexes on a
curve ($p(F_1(R(U)))$)
 such that the set of generators
of the cones ($p(R(U))$) has the
corresponding vertex as
generic focal locus in each line.}\\ 
\cline{2-3}
{ \quad}& {\footnotesize $\alpha_3$}& {\footnotesize $1$--dimensional
families of planes containing each of the
generators of an irreducible ruled surface ($p(F_1(R(U)))$), 
 nondevelopable, except the family of tangent planes to
such a surface.}\\
\hline { \quad}& {\footnotesize $\beta_1$}& {\footnotesize Family of tangent
planes to an irreducible nondegenerate surface
without asymptotic curves ($\Sigma'\subset{\bf P}^4$).}\\
\cline{2-3}
 {\footnotesize $\beta$: Two different developable systems passing through the
general plane. }&
 {\footnotesize $\beta_2$}& {\footnotesize A
$1$--dimensional family of planes passing through each
one of the tangent lines of a curve ($\Sigma'$).
If these families are linear they cannot contain
the corresponding osculating plane to
$\Sigma'$.}\\ \cline{2-3}
{ }& {\footnotesize $\beta_3$}& {\footnotesize Cone (with vertex a
point $\Sigma'$) of a nonparabolic congruence
of lines in a hyperplane $H\subset{\bf P}^4$,
$\Sigma'\notin H$.}\\ 
\hline
{ \quad}& {\footnotesize $\gamma_1$}& {\footnotesize Family of tangent planes to an
irreducible nondegenerate surface
with asymptotic curves ($\Sigma'\subset{\bf P}^4$).}\\
\cline{2-3}
 {\footnotesize $\gamma$: Two coincident developable systems passing through
the general plane.}&
 {\footnotesize $\gamma_2$}&{\footnotesize A
$1$--dimensional linear family of planes passing 
through each one of the tangent lines to a curve
($\Sigma'$) satisfying that each one of such
families must contain the corresponding osculating
plane to $\Sigma'$.}\\ \cline{2-3}
{ \quad}& {\footnotesize $\gamma_3$}& {\footnotesize Cone (with vertex a
point $\Sigma'$) of a parabolic congruence of
lines in a hyperplane $H\subset{\bf P}^4$,
$\Sigma'\notin H$.}\\ \hline  {\footnotesize $\delta$: An infinity of developable systems 
passing through the general plane.}& {\footnotesize $\delta$}& {\footnotesize
Family of the planes in
${\bf  P}^4$ containing a line ($p(F_1(\Sigma))$).}\\
\hline 
\end{tabular}

\end{document}